\documentclass{article}

\usepackage{arxiv}

\usepackage[utf8]{inputenc} 
\usepackage[T1]{fontenc}    
\usepackage{hyperref}       
\usepackage{url}            
\usepackage{booktabs}       
\usepackage{amsfonts}       
\usepackage{nicefrac}       
\usepackage{microtype}      
\usepackage{lipsum}		    
\usepackage{graphicx}
\usepackage[square,sort,comma,numbers]{natbib}
\usepackage{doi}
\usepackage{float}
\usepackage{caption}
\usepackage{subcaption}
\usepackage{xcolor}
\usepackage{amsmath}
\usepackage{svg}
\usepackage{enumitem}

\newcommand{\define}[4]{$#1=#1(#2):#3\mapsto#4$}

\newcommand{\IPeq}[3]{\big \langle#1,#2 \big \rangle_{#3}}
\newcommand{\R}{\mathbb{R}}

\title{FNPF-SEM: A parallel spectral element model in Firedrake for fully nonlinear water wave simulations}

\author{ \href{https://orcid.org/0000-0001-6698-2623}{\includegraphics[scale=0.06]{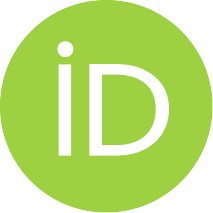}\hspace{1mm}Jens Visbech} \\
	Department of Applied Mathematics and Computer Science\\
	Technical University of Denmark\\
	Kongens Lyngby, 2800, Denmark \\
	\texttt{jvis@dtu.dk} \\
\And
	\href{https://orcid.org/0000-0003-3310-9638}{\includegraphics[scale=0.06]{Figures_pdf/orcid-eps-converted-to.pdf}\hspace{1mm}Anders Melander} \\
    Department of Applied Mathematics and Computer Science\\
	Technical University of Denmark\\
	Kongens Lyngby, 2800, Denmark \\
	\texttt{adame@dtu.dk} \\
\And
    	\href{https://orcid.org/0000-0001-8626-1575}{\includegraphics[scale=0.06]{Figures_pdf/orcid-eps-converted-to.pdf}\hspace{1mm}Allan Peter Engsig-Karup} \\
	Department of Applied Mathematics and Computer Science\\
	Technical University of Denmark\\
	Kongens Lyngby, 2800, Denmark \\
	\texttt{apek@dtu.dk} \\
}

\renewcommand{\headeright}{Visbech et al. (2025)}
\renewcommand{\shorttitle}{An \textit{arXiv} Preprint}

\hypersetup{
pdftitle={},
pdfsubject={},
pdfauthor={},
pdfkeywords={},
}

\begin{document}
\maketitle

\begin{abstract}

We present a new parallel spectral element solver, FNPF-SEM, for simulating linear and fully nonlinear potential flow-based water waves and their interaction with offshore structures. The tool is designed as a general-purpose wave model for offshore engineering applications.
Built within the open-source framework Firedrake, the new FNPF-SEM model is designed as a computational tool capable of capturing both linear and nonlinear wave phenomena with high accuracy and efficiency, with support for high-order (spectral) finite elements. Additionally, Firedrake provides native support for MPI-based parallelism, allowing for efficient multi-CPU distributed computations needed for large-scale simulations.
We demonstrate the capabilities of the high-order spectral element model through $h$- and $p$-convergence studies, and weak and strong scaling tests. Validation is performed against analytical solutions and experimental data for several benchmark cases, including nonlinear high-order harmonic generation and linear and nonlinear wave interactions with a cylinder and a breakwater.
The new FNPF-SEM model offers a numerical framework for simulating wave propagation and wave–structure interactions, with the following key features: i) the ability to represent complex geometries through flexible, unstructured finite element meshes; ii) reduced numerical diffusion and dispersion by using high-order polynomial expansions; and iii) scalability to full- and large-scale simulations over long time periods through a parallel implementation.
\end{abstract}

\keywords{Spectral element method \and linear and nonlinear free surface potential flow \and water waves \and wave-structure interactions \and wave propagation \and parallel computing}

\section{Introduction}\label{sec:intro}

Simulating the dynamics of water waves remains a significant challenge in computational fluid dynamics (CFD), especially when considering complex geometry due to the flow domain and included marine or offshore structures. However, despite the challenge, accurate modeling of wave propagation and interactions with bathymetry and structures is essential for marine and offshore engineering applications. For example, estimating how the local wave climate and the resulting wave-induced forces impact offshore structures such as large offshore wind turbines  \cite{bailey2015metocean}, wave energy devices, and breakwaters. Other cases include metocean studies such as assessments of ocean conditions, \cite{bitnergregersen2015joint}. Many offshore engineering applications have to consider modeling of nonlinear water waves in large regions and within large time scales, requiring effective parallel implementations of such wave models to accurately and efficiently compute real-life scenarios within acceptable computational times.

Physical accuracy is typically associated with the chosen mathematical model and how it represents the physics. High-fidelity simulations can be based on the Navier-Stokes equations as a model for all fluid flows. However, when solved numerically, these equations come with high computational requirements. For offshore engineering applications, a broad range of free surface flow applications can be approximated by fully nonlinear potential flow (FNPF), assuming the fluid is incompressible, inviscid, and irrotational. Moreover, if the wave steepness is small, the linear potential flow (LPF) equations are applicable, i.e., a small amplitude approximation, which reduces the problem complexity. One desirable feature of the FNPF formulation is the ability to simulate both offshore and nearshore environments, since wave interactions with the bathymetry are accurately represented. When the simulated problem satisfies the assumptions of FNPF theory, a performance improvement of at least two orders of magnitude can be achieved while achieving an accuracy comparable to the full CFD simulations, as shown in \cite{ransley2019blind}. If the problem adheres further to the assumptions of LPF theory, additional orders of magnitude in computational efficiency can be realized.

Numerically accurate and efficient solutions to the FNPF and LPF formulations have been a major research topic for decades. Traditionally, the boundary integral method has been popular for wave modeling. To highlight a few of the well-known and recent boundary element-based models (including the traditional low-order boundary element method, the quadratic boundary element method, and the high-order boundary element method), see \cite{liu2001computations, xue2001computations,lee2006wamit,yan2011efficient,harris2022nonlinear,seixasdemedeiros2024fast,kurnia2023nemoh}. Another option are high-order spectral methods, which are very efficient when modeling nonlinear wave propagation and submerged wave-structure interactions, \cite{dommermuth1987high,liu1992high,christiansen2013efficient,melander2024multigrid_2}. The FNPF and LPF equations have also been solved using finite difference schemes, \cite{bingham2007accuracy,engsigkarup2009efficient,ducrozet2010high,amini2017solving,amini2018pseudo}, and using harmonic polynomials, \cite{shao2014harmonic,tong2024adaptive}. Using the classical second-order accurate finite element method, we have the work by \cite{wu1994finite} and many more.

Spectral element methods (SEMs), initially introduced by \cite{patera1984spectral}, combine key advantages of both spectral (polynomial) and finite element methods. SEMs allow for two distinct convergence strategies: (i) $h$-convergence, where the mesh is refined while maintaining a fixed polynomial order, resulting in algebraic convergence; and (ii) $p$-convergence, where the mesh is held fixed while increasing the polynomial order of the basis functions, typically leading to spectral (exponential) convergence. The hybrid approach provides efficient and highly accurate numerical solutions due to the use of piecewise polynomial basis functions of arbitrary order with spectral accuracy. Moreover, the SEM can handle complex geometries common in engineering applications because it uses unstructured mesh tessellations. This combination of high accuracy and geometrical flexibility addresses traditional numerical approaches' limitations, e.g., difficulty in modeling complex geometry and/or high dispersive and diffusive errors. For further details and applications of the SEM, see the review by \cite{xu2018spectral} and the textbook by \cite{karniadakis2005spectral}. 

Many open-source finite element frameworks that allow for high-order basis functions and parallel computing have been developed over the years. For example, the spectral/$hp$ element framework Nektar++, \cite{cantwell2015nektar}, the arbitrary high-order finite element framework MFEM, \cite{anderson2021mfem}, the deal.II finite element framework, \cite{arndt2021dealii}, or libParanumal, \cite{chalmers2020libParanumal}.
Another computational framework that allows for high-order basis functions is Firedrake \cite{ham2023firedrake}. Firedrake is an open-source framework that leverages itself around automating the required steps to solve partial differential equations using finite elements. One of Firedrake's main features is the unified form language, \cite{alnaes2014unified}, which allows for the symbolic expression of, e.g., any weak variational formulations. Moreover, Firedrake uses state-of-the-art PETSc solvers \cite{balay1997efficient,balay2024usersmanual}, and supports structured and unstructured meshes using the most common types of elements and finite element spaces. Lastly, Firedrake allows for a fully parallelized computational setup through MPI and can run on large many-core hardware systems. In the recent decade, the Firedrake project has established itself as a robust open-source platform for implementing advanced numerical models.

At DTU, SEM-based FNPF models for simulating free surface water waves have been under active development for more than a decade. This effort began with the work of  \cite{engsigkarup2016stabilised}, who proposed a solution to the intrinsic mesh-induced instability associated with finite element discretizations of free surface flows, as previously identified in \cite{robertson1999free}. The SEM implementations at DTU are based primarily on the FNPF formulations, \cite{engsigkarup2016stabilised, engsigkarup2019spectral, engsigkarup2019mixed, engsigkarup2021efficient,visbech2024high_2}, and the LPF approximations, \cite{visbech2024high_1,visbech2024solving}. This combined FNPF-SEM approach offers a numerical approach that can handle complex geometries and simulate over long time periods accurately by using the high-order polynomial basis. 
The methodology has been extended to encompass the full incompressible Navier–Stokes equations in \cite{melander2024multigrid_1} with a free surface. Achieving low dispersive and diffusive numerical errors that stem from the numerical discretization procedure is a key aspect for long-time wave simulations. High-order discretization methods, such as SEM, are very efficient compared to low-order methods \cite{kreiss1972comparison,engsigkarup2016stabilised}.
In addition, \cite{eskilsson2024water} presented a high-performance implementation of an FNPF-based water wave solver using Nektar++. The current capabilities of these FNPF-based SEM models are comparable in functionality, benefiting significantly from the robust and scalable features provided by modern open-source computational frameworks. Moreover, efficient water wave solutions to the Benney-Luke-type equations stated using variational principles are also available via Firedrake, \cite{bokhove2016variational}. Furthermore, the Thetis project, an unstructured grid coastal ocean model, was built using the Firedrake framework, \cite{karna2018thetis}.

\subsection{Paper contribution}

The main aim is to verify and validate a new parallel implementation of the aforementioned DTU-based SEM models for solving FNPF and LPF wave problems. The parallel code is implemented in  Firedrake, allowing for an efficient parallel implementation of the FNPF equations to support large- and full-scale wave problems over long time horizons, ultimately to be used within offshore engineering applications.
The paper verifies the numerical implementation, referred to as FNPF-SEM, through $h$- and $p$-convergence studies and the parallel setup through strong and weak scaling tests. Moreover, the entire model is validated in various test cases for FNPF and LPF wave propagation and wave structures interactions.

\subsection{Outline of paper}

Following the introduction in Section \ref{sec:intro}, the paper proceeds with the mathematical model in Section \ref{sec:math_model}, which describes the computational domain, its boundaries, and the governing equations. Subsequently, Section \ref{sec:num_dis} outlines the numerical discretization approaches, highlighting the chosen spatial and temporal discretization approaches, including key aspects such as wave initialization, generation, and absorption. The results are shown in Section \ref{sec:res}, which presents both numerical verification studies (including convergence and scaling analyses) and validation cases involving linear and nonlinear wave propagation and wave–structure interactions. Finally, the paper concludes with a summary of the main findings in Section \ref{sec:con}. An appendix is also included for supplementary material in Section \ref{sec:app}.
\section{Mathematical Model}\label{sec:math_model}

\begin{figure}[H]
    \centering
    \includegraphics[scale=0.95]{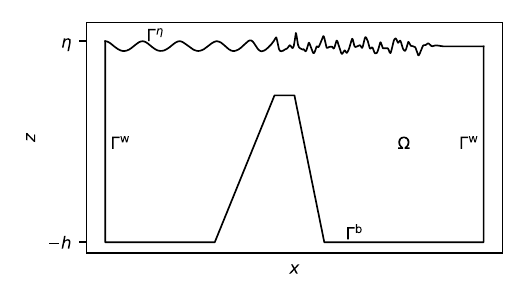}
    \caption{Conceptual schematic of the fluid domain, $\Omega$, and its boundaries in the $xz$-plane.}
    \label{fig:FluidDomain}
\end{figure}

Considering a fluid domain, $\Omega \in \mathbb{R}^d$, where $d = 3$, and its completely surrounding boundary, $\Gamma = \partial \Omega$. Here, $\boldsymbol{x} = (x,y,z)$ where $x$ and $y$ span the horizontal plane, and $z$ is the vertical direction defined positively upwards. For the boundary, we have the following decomposition: $\Gamma = \bigcup_i \Gamma^{i}$ for $i \in \{\eta, \text{b}, \text{w},\text{p} \}$. Here, $\Gamma^{\eta}$ is the top free surface boundary defined at $z = \eta$ in the FNPF setting, where \define{\eta}{x,y,t}{\Gamma^{\eta} \times \mathcal{T}}{\R} is the free surface elevation defined from $z = 0$ in the time-domain $\mathcal{T}:t \geq 0$. In the LPF setting, $\Gamma^{\eta}$ is defined at $z = 0$. Moreover, $\Gamma^{\text{b}}$ is the bathymetry boundary at the bottom at $z = -h$, where \define{h}{x,y}{\Gamma^{\text{b}}}{\R} is the height of the water column measured around $z = 0$. At the vertical boundaries, two possibilities exist: $\Gamma^{\text{w}}$ denotes impermeable walls, and $\Gamma^{\text{p}}$ denotes periodic boundaries. A conceptual sketch of $\Omega$ and its boundaries is shown in Figure \ref{fig:FluidDomain} for the $xz$-plane. \textbf{Remark:} Partially or fully submerged structures can be included in this setting, yet this is not the scope of this work.

\subsection{Governing equations}

Assuming potential flow, i.e., an incompressible, irrotational, and inviscid fluid, the gradient of a scalar velocity potential, \define{\phi}{\boldsymbol{x},t}{\Omega \times \mathcal{T}}{\R}, can be used to express the fluid velocities, $\boldsymbol{u} = (u,v,w) = (u(\boldsymbol{x},t),v(\boldsymbol{x},t),w(\boldsymbol{x},t)) : \Omega \times \mathcal{T} \mapsto \R^d$ as $\boldsymbol{u} = \boldsymbol{\nabla} \phi$ where $\boldsymbol{\nabla} = (\partial_x, \partial_y, \partial_z)$ is the Cartesian gradient operator in $\Omega$. The spatial mathematical problem is to find $\phi \in C^2(\Omega)$ by solving the Laplace equation with suitable boundary conditions as
\begin{align}
\begin{split}\label{eq:Laplace_strong}
     \boldsymbol{\nabla}^2 \phi  &= 0, \quad \text{in} \quad \Omega, \\
        \phi &= \phi_{\eta}, \quad \text{on} \quad \Gamma^{\eta}, \\
      \boldsymbol{\nabla} \phi \cdot \boldsymbol{n} &= 0, \quad \text{on} \quad \{ \Gamma^{\text{b}}, \Gamma^{\text{w}} \}, \\
     \phi |_{\Gamma^{\text{p}}} &= \phi |_{\Gamma^{\text{p}}}, \quad \text{on} \quad \Gamma^{\text{p}},
     \end{split}
\end{align}
where $\boldsymbol{n} = (n_x,n_y,n_z) = (n_x(\boldsymbol{x}),n_y(\boldsymbol{x}),n_z(\boldsymbol{x})): \Gamma \mapsto \R^d$, is the outward-facing normal vector to any boundary of interest. The subscript notation of "$\eta$" denotes quantities defined on $\Gamma^{\eta}$, i.e., \define{\phi_{\eta}}{x,y,t}{\Gamma^{\eta} \times \mathcal{T}}{\R}. This notation becomes useful when we consider the FNPF time governing kinematic and dynamic free surface conditions on Zakharov form, \cite{zakharov1968stability}. On $\Gamma^{\eta} \times \mathcal{T}$, we have 
\begin{align}
    \partial_t \eta &= - \boldsymbol{\widetilde{\nabla}} \eta \cdot \boldsymbol{\widetilde{\nabla}} \phi_{\eta} + w_{\eta} (1 + \boldsymbol{\widetilde{\nabla}} \eta \cdot \boldsymbol{\widetilde{\nabla}} \eta), \label{eq:kinematic_free_surface_strong} \\
    \partial_t \phi_{\eta} &= - g \eta - \frac{1}{2} \boldsymbol{\widetilde{\nabla}} \phi_{\eta} \cdot \boldsymbol{\widetilde{\nabla}} \phi_{\eta} +  \frac{1}{2} w_{\eta}^2 (1 +  \boldsymbol{\widetilde{\nabla}} \eta \cdot \boldsymbol{\widetilde{\nabla}} \eta). \label{eq:dynamic_free_surface_strong}
\end{align}
Here $\boldsymbol{\widetilde{\nabla}} = (\partial_x,\partial_y)$ is the Cartesian gradient operator only defined on the horizontal free surface, $\Gamma^{\eta}$, and $g = 9.82$ [m/s$^2$] is the gravitational acceleration. Given initial conditions for $\eta$ and $\phi_{\eta}$ at $t=0$, this coupled FNPF initial boundary value problem can be solved uniquely for any point in time and space. \textbf{Remark:} The classical linearized version of \eqref{eq:Laplace_strong} through \eqref{eq:dynamic_free_surface_strong} is considered as well; however, these governing equations are omitted for brevity.
\section{Numerical Discretization}\label{sec:num_dis}

We adopt a classical \textit{method of lines} approach to the governing equations. The spatial and temporal discretizations are discussed below individually. As stated earlier, the mathematical model is implemented in Firedrake, \cite{ham2023firedrake}. Firedrake is an open-source computational framework for solving partial differential equations using the finite element method. It provides a high-level interface based on the unified form language, \cite{alnaes2014unified}, allowing users to express weak variational formulations in a compact and readable form. At the backend, Firedrake automates the generation of optimized C code and supports efficient parallel execution via MPI, making it well-suited for large-scale simulations. Built on PETSc, \cite{balay2024usersmanual}, it includes state-of-the-art solvers and supports a wide range of element types and mesh structures, including high-order spectral elements on both structured and unstructured mesh tessellations.

\subsection{Spatial discretization}

We perform a discrete tesselation of the spatial domain, $\Omega$, and boundary, $\Gamma^{\eta}$, using non-overlapping finite elements as 
\begin{equation*}
    \mathcal{D} \simeq \mathcal{T}_{\mathcal{D}} = \bigcup_{n=1}^{N_{\mathcal{D}}} \mathcal{D}_n, \quad \text{for} \quad \mathcal{D} \in \{\Omega, \Gamma^{\eta} \},
\end{equation*}
where $N_{\mathcal{D}}$ is an arbitrary number of elements, such that $\mathcal{D}_n$ is the $n$'th element in the tesselation $\mathcal{T}_{\mathcal{D}}$. The element types and meshes are highlighted in the following section. On the two tesselations, we introduce $C^0(\mathcal{D})$ function spaces of piecewise polynomial functions as:
\begin{equation*}
    \mathcal{V}_\mathcal{D} = \{\nu_{\mathcal{D}} \in C^0(\mathcal{D}); \forall n \in \{ 1,..., N_{\mathcal{D}} \}, \nu_{\mathcal{D}}|_{\mathcal{D}_n} \in \mathbb{P}^p\},
\end{equation*}
where $\mathbb{P}^p$ denotes the polynomial space of nodal Lagrange polynomials of degree at most $p$.

\subsubsection{Weak variational formulations} 

Next, we state the weak forms of \eqref{eq:Laplace_strong} and \eqref{eq:kinematic_free_surface_strong} through \eqref{eq:dynamic_free_surface_strong}, respectively. The linear variational problem for the Laplace problem is to find $\phi \in \mathcal{V}_{\Omega}$ such that $\forall \nu_{\Omega} \in \mathcal{V}_{\Omega}$, with
\begin{equation}\label{eq:Laplace_weak}
    a[\phi;\nu_{\Omega}] = b[\nu_{\Omega}],
\end{equation}
where the bi-linear and linear forms are derived by multiplying by a test function, $\nu_{\Omega}$, integrating over $\Omega$, and applying integration by parts as
\begin{align*}
    a[\phi;\nu_{\Omega}] &= \IPeq{\boldsymbol{\nabla} \phi}{\boldsymbol{\nabla} \nu_{\Omega}}{\Omega}, \\
    b[\nu_{\Omega}] &= 0,
\end{align*}
with $\langle a,b \rangle_S = \int_S a \cdot b ~ d S$ denoting the inner product over the set $S$.

For the two free surface conditions, we define a nonlinear variational problem to find $\{\partial_t \eta, \partial_t  \phi_{\eta} \} \in \mathcal{V}_{\Gamma^{\eta}}$ such that $\forall \nu_{\Gamma^{\eta}} \in \mathcal{V}_{\Gamma^{\eta}}$, we solve
\begin{equation}\label{eq:free_surface_weak}
    F_k[\partial_t \eta; \nu_{\Gamma^{\eta}}] = 0, \quad \text{and} \quad F_d[\partial_t \phi_{\eta}; \nu_{\Gamma^{\eta}}] = 0,
\end{equation}
where the nonlinear variational forms are derived by multiplying by a test function, $\nu_{\Gamma^{\eta}}$, integrating over $\Gamma^{\eta}$, and moving everything to the left-hand side. Ultimately, this gives
\begin{equation*}
    \begin{split}
        &F_k[\partial_t \eta; \nu_{\Gamma^{\eta}}] = \IPeq{\partial_t \eta}{\nu_{\Gamma^{\eta}}}{\Gamma^{\eta}} +         \IPeq{\boldsymbol{\widetilde{\nabla}} \eta \cdot \boldsymbol{\widetilde{\nabla}} \phi_{\eta} - w_{\eta} (1 + \boldsymbol{\widetilde{\nabla}} \eta \cdot \boldsymbol{\widetilde{\nabla}} \eta)}{\nu_{\Gamma^{\eta}}}{\Gamma^{\eta}},
    \end{split}
\end{equation*}
and
\begin{equation*}
    \begin{split}
        &F_d[\partial_t \phi_{\eta}; \nu_{\Gamma^{\eta}}] = \IPeq{\partial_t \phi_{\eta}}{\nu_{\Gamma^{\eta}}}{\Gamma^{\eta}} + 
        \IPeq{g \eta - \frac{1}{2} \boldsymbol{\widetilde{\nabla}} \phi_{\eta} \cdot \boldsymbol{\widetilde{\nabla}} \phi_{\eta} -  \frac{1}{2} w_{\eta}^2 (1 +  \boldsymbol{\widetilde{\nabla}} \eta \cdot \boldsymbol{\widetilde{\nabla}} \eta)}{\nu_{\Gamma^{\eta}}}{\Gamma^{\eta}}.
  \end{split}
\end{equation*}

\subsubsection{Initial meshes, extrusion, and mesh updating}

\begin{figure}[H]
    \centering
    \includegraphics[scale=0.95]{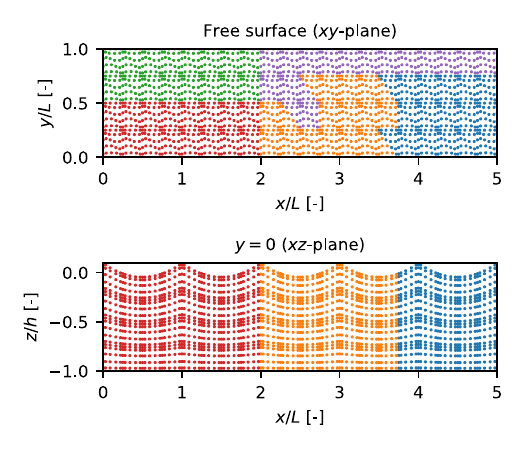}
    \caption{Point distribution (indicated by different colors) for a finite FNPF wave tank distributed over $N_c = 5$ cores on a mesh with $(N_x,N_y,N_z) = (20,4,4)$ elements of order $p=5$.}
    \label{fig:CoordinatePlot}
\end{figure}

Firedrake has native support for extruded meshes, \cite{mcrae2016automated, bercea2016structure}. The mesh construction for this model proceeds in three steps: 
First, a triangular mesh for the free surface, $\mathcal{T}_{\Gamma^{\eta}}$, is generated. This is done using Firedrake’s built-in mesh utilities for simple rectangular domains or via Gmsh for more complex geometries, \cite{geuzaine2009gmsh}.
Second, the surface mesh is vertically extruded in the $z$-direction into $N_z$ layers to form the three-dimensional fluid domain mesh, $\mathcal{T}_{\Omega}$. These two meshes, referred to as the "initial meshes" $\big(\mathcal{T}_{\Gamma^{\eta}}^i$ and $\mathcal{T}_{\Omega}^i \big)$, consist entirely of prismatic elements, and can be unstructured to adapt the geometric features in the horizontal $(x,y)$-plane and structured in the $z$-direction. 
Third, in the FNPF setting, a vertical translation is applied to the mesh along each column of points to track the evolving free surface curvature. This ensures an accurate geometric representation of the time-dependent fluid domain, $\Omega$, as defined by the free surface elevation, $\eta$. Consequently, mesh updates are performed at each discrete time stage or step to reflect the dynamic boundary. In Figure \ref{fig:CoordinatePlot} (lower), the updated curvilinear mesh can be seen. \textbf{Remark:} A similar approach is applied to modify the mesh due to non-flat bathymetries. Currently, submerged and floating structures are not directly supported in the model due to the mesh extrusion; however, work is being done to include this.

We want to move information and data (when running in parallel) between the free surface mesh, $\mathcal{T}_{\Gamma^{\eta}}$, and the corresponding top free surface points of the fluid domain mesh, $\mathcal{T}_{\Omega}$. Using extruded meshes in Firedrake, it is guaranteed that the parallelized partition of the two meshes coincides (in the $x$- and $y$-coordinates) as seen in Figure \ref{fig:CoordinatePlot}. This feature makes it possible to compute a mapping, $m: \mathcal{T}_{\Omega} \mapsto \mathcal{T}_{\Gamma^{\eta}}$, that extracts the free surface values, $\boldsymbol{f_{\Gamma^{\eta}}}$, from any fluid domain variable, $\boldsymbol{f_{\Omega}}$, as $\boldsymbol{f_{\Gamma^{\eta}}} = \boldsymbol{f_{\Omega}}[m]$, i.e., a simple index/mask operation between two vectors.

\subsection{Temporal discretization and modal filtering}

Time integration is performed using a standard explicit 4th-order, 4-stage, Runge-Kutta scheme (ERK4). The time step size, $\Delta t$, is kept constant throughout the simulation and is computed based on a Courant-Friedrichs-Lewy (CFL) condition: $\Delta t = C_{\text{CFL}} \Delta x_{\min} u_{\max}^{-1}$. Here, $C_{\text{CFL}} = 0.95$, $\Delta x_{\min}$ denotes the minimum spacing between two grid points on $\Gamma^{\eta}$, and $u_{\max}$ represents the maximum velocity encountered in the simulation.

Some form of stabilization is typically required to model strongly nonlinear wave phenomena accurately. In \cite{engsigkarup2016stabilised}, the instability issue initially highlighted by \cite{robertson1999free} was addressed using structured elements aligned in the vertical direction. This was coupled with an over-integration (collocation) approach for the nonlinear free surface terms, \cite{kirby2003dealiasing}, and a mild modal (spectral) filter applied to the free surface quantities, \cite{hesthaven2008filtering}. In Firedrake, the integration order of individual integral terms can be explicitly specified, facilitating the over-integration strategy suitable for reducing aliasing instabilities of nonlinear wave propagation, \cite{engsigkarup2016stabilised}. A modal filter was implemented within this framework by developing a custom \texttt{PyOP2} kernel written in C, executed via the \texttt{par\_loop} interface. This allows for the parallel application of the filter across all free surface elements on the interior point values. Specifically, the kernel operates element-wise, applying the local modal filter as defined in \cite{engsigkarup2016stabilised}: 
Suppose $g_n(x)$ is a function defined on a free surface element $\Gamma^{\eta}_n$. On each element, we have $n_p = (p+1)(p+2)/2$ nodal points denoted by $\boldsymbol{x_n} = (x_n,y_n)$. We map the solution to the modal space as
\begin{equation*}
    \widehat{g}_n(f) = \mathcal{V}^{-1} ~ g_n(\boldsymbol{x_n}),
\end{equation*}
where $\mathcal{V} \in \mathbb{R}^{n_p \times n_p}$ is a conforming Vandermonde matrix defined as $\mathcal{V}_{ij} = \widehat{\psi}_j(\boldsymbol{r_i})$, where $\boldsymbol{r_i} = (r_i,s_i)$, is coordinates in a reference domain, $\mathcal{R} = \{\boldsymbol{r} = (r,s) | (r,s) \geq -1; r+s \leq 0 \}$ and $\widehat{\psi}(\boldsymbol{r})$ is a modal basis function. See \cite{engsigkarup2016stabilised} for more details. Using a filter matrix, $\mathcal{F} \in \mathbb{R}^{n_p \times n_p}$, filtering can be applied to the highest modes
\begin{equation*}
    \widehat{g}_{n,filtered}(f) = \mathcal{F} ~ \widehat{g}_n(f).
\end{equation*}
Then, the filtered solution can be mapped back to nodal space
\begin{equation*}
    g_{n,filtered}(\boldsymbol{x_n}) = \mathcal{V} ~ \widehat{g}_{n,filtered}(f).
\end{equation*}
This filtering approach can be applied to any free surface quantity, e.g., $\eta$, $\phi_{\eta}$, or $w_{\eta}$, as needed for nonlinear waves.

\subsection{Wave initialization, generation, and absorption}

For periodic wave tank simulations, the model is initialized with $\eta$ and $\phi_{\eta}$, either from the linear wave theory in the LPF setting or from an accurate stream function solution in the FNPF setup, \cite{fenton1988numerical}. Linear (small-amplitude) waves can be fully characterized by the non-dimensional wave number, $kh$, given a specified wavelength, $L$, or water depth, $h$. In contrast, nonlinear (finite-amplitude) stream function waves also require the specification of a relative wave steepness, $\varepsilon / \varepsilon_{\max}$, to quantify their nonlinearity. Here, $\varepsilon = H/L$ denotes the wave steepness based on the wave height $H$, and $\varepsilon_{\max}$ is the maximum steepness attainable before wave breaking for a given wave, \cite{fenton1990nonlinear}.
The initial conditions are set to zero for simulations in numerical wave tanks of finite size. Wave generation and absorption are handled using classical relaxation zones, \cite{larsen1983open}, of 2-3 wavelengths in size. The wave fields, $\eta$ and $\phi_{\eta}$, are then modified at each time step as a post-processing operation as
\begin{equation*}
    f = (1-C_r) f + C_r f_{\text{true}},
\end{equation*}
where $f \in \{\eta, \phi_{\eta} \}$ and \define{C_r}{x,y}{\Gamma^{\eta}}{\R} is the relaxation function, taking the value of zero in the target zone and smoothly goes to unity in the generation and absorption zone. In the generation zone, $f_{\text{true}}$ is the true actual solution and zero outside this zone. To avoid impulsive behavior at the onset of the simulation, the generation zone is gradually activated from 0\% to 100\% over a duration of five wave periods.

\subsection{Entire computational workflow}

To support the understanding of how the entire computational workflow is outlined, a \textit{step-by-step} procedure is shown below:

\begin{itemize}[leftmargin=*]

    \item[] \textbf{Step 1 -- Setup:} Build initial meshes, $\mathcal{T}_{\Gamma^{\eta}}^i$ and $\mathcal{T}_{\Omega}^i$, and compute pre-computable quantities, e.g., the map, $m$.
    
    \item[] \textbf{Step 2 -- Initial conditions:} Impose initial conditions on $\eta$ and $\phi_{\eta}$ at $t = t_{\text{start}}$.

    \item[] \textbf{Step 3 -- Time integration:} Iterate \textbf{Step 4} below over time steps/stages until $t = t_{\text{end}}$.

    \begin{itemize}[leftmargin=*]

    \item[] \textbf{Step 4.1 -- Mesh-update:} Update $\mathcal{T}_{\Omega}$ according to $\eta$.

    \item[] \textbf{Step 4.2 -- Boundary conditions:} Impose direct Dirichlet boundary conditions using $\phi_{\eta}$.

    \item[] \textbf{Step 4.3 -- Laplace problem:} Solve the Laplace problem for the velocity potential, $\phi$.

    \item[] \textbf{Step 4.4 -- Gradient recovery:} Compute the vertical velocity, $w$, from $\phi$ via gradient recovery.

    \item[] \textbf{Step 4.5 -- Vertical free surface velocity:} Extract the vertical free surface velocity, $w_{\eta}$, from $w$ via $m$.

    \item[] \textbf{Step 4.6 -- Free surface conditions:} Evaluate the two free surface conditions using $\eta$, $\phi_{\eta}$, and $w_{\eta}$.

    \item[] \textbf{Step 4.7 -- Relaxation and filtering:} Modify wave solution depending on generation and absorption zones and apply filtering as needed.

    \end{itemize}

    \item[] \textbf{Step 5 - Post-processing:} Plot results, etc.
    
\end{itemize}

\section{Results}\label{sec:res}

In the following, we present numerical results related to the verification and validation of the model. All simulations are performed on the HPC system at the DTU Computing Center at the Technical University of Denmark using the hardware configurations listed in Table \ref{tab:Configurations}. 

\begin{table}[H]
    \centering
    \caption{Hardware configurations.}
    \begin{tabular}{l l} \hline
        Total nodes, $N_n$, available & 15 \\ 
        Total cores, $N_c$, available & 300 \\
        Cores per node, $N_c/N_n$ & 20 \\ 
        Memory per node & 125 GB \\
        CPU & Intel Xeon E5-2660 v3 \\
        CPU frequency & 2.60 GHz \\ 
        \hline
    \end{tabular}
    \label{tab:Configurations}
\end{table}

As for the numerical solver strategies, we solve the Laplace problem in \eqref{eq:Laplace_strong} following \cite{rathgeber2016firedrake} by using a CG solver with the HYPRE BoomerAMG algebraic multigrid preconditioner, \cite{hestenes1952methods,falgout2006design}, using a relative and absolute tolerance of $10^{-6}$ and $10^{-15}$, respectively. For the evaluation of the two nonlinear free surface conditions in \eqref{eq:kinematic_free_surface_strong} and \eqref{eq:dynamic_free_surface_strong}, we employ the default solver using a GMRES solver with a ILU as preconditioner and a relative and absolute tolerance of $10^{-5}$ and $10^{-15}$, respectively. For the purpose of this study, this leads to an efficient solver strategy for the cases considered, cf. the analysis of iterative solver efficiency for FNPF modeling reported in \cite{engsigkarup2014analysis}.

\subsection{Model verification}

In the verification of the model, various timings will be presented. Those timings are an average for the number of cores, $N_c$, used in the test. When using $N_c \leq 20$, a complete node is allocated to the run. 

\subsubsection{Convergence \textit{hp} studies}

Investigating and confirming the expected convergence rates of the proposed model constitutes the first step in the verification process. For this, we perform $h$- and $p$-convergence studies; that is, mesh refinement under a fixed polynomial order of the basis functions, and an increase in the polynomial order on fixed meshes, respectively. For the analysis, a periodic domain is considered, featuring nonlinear waves for all combinations of $kh = \{1, 3, 6\}$ and $\varepsilon / \varepsilon_{\max} = \{10\%, 50\%, 90\% \}$. The error is evaluated in the $\infty$-norm of the vertical free surface velocity, $|w_{\eta}^{\text{true}} - w_{\eta}|_{\infty}$, where $w_{\eta}^{\text{true}}$ denotes the reference (true) solution and $w_{\eta}$ is the corresponding numerical solution. All wave solutions are defined over a domain of fixed length $L = 1$ [m], implying $L_x = L_y = L$. The depth is varied depending on the specific value of $kh$.

\begin{figure}[H]
    \centering
    \includegraphics[scale=0.90]{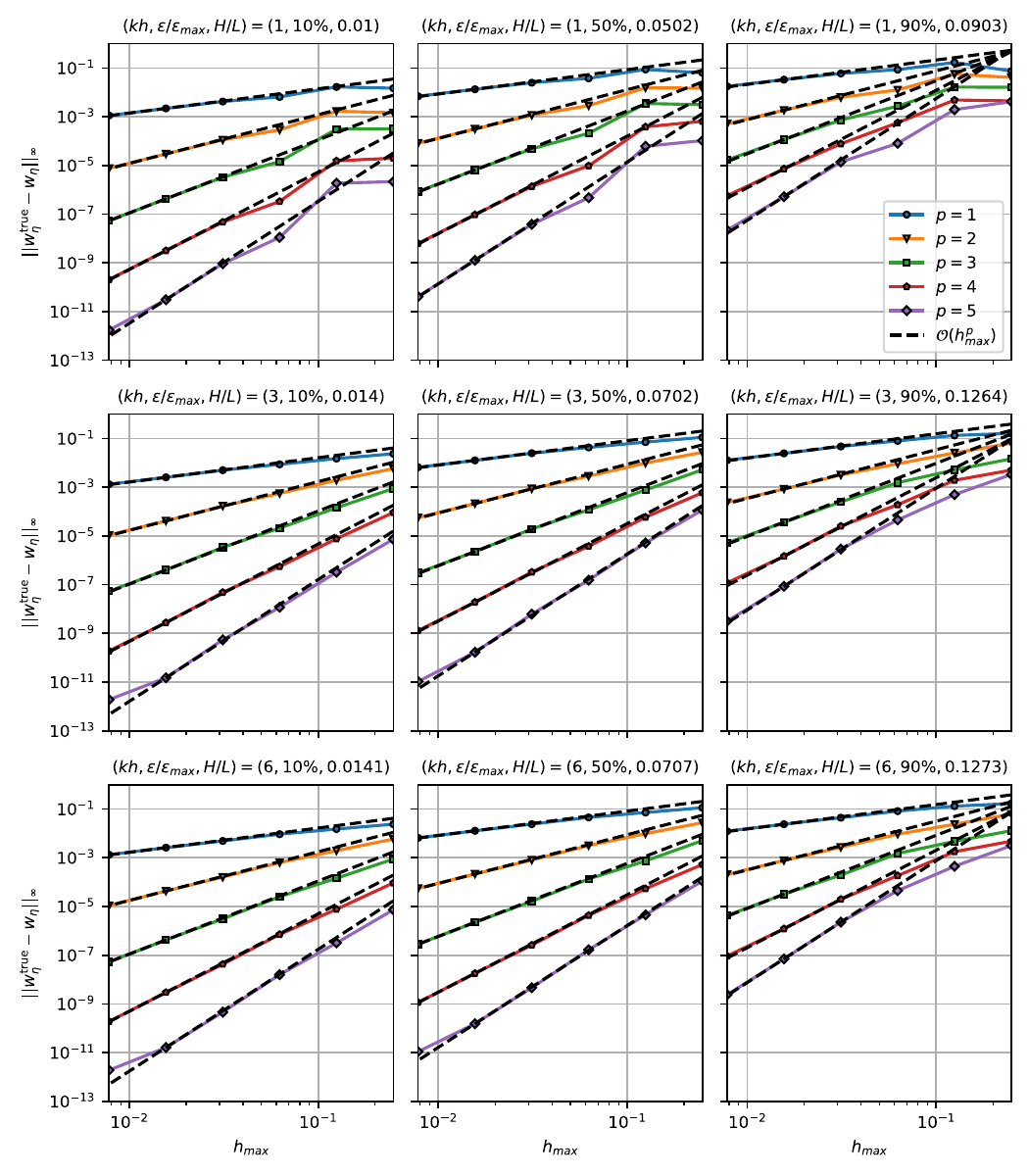}
    \caption{Convergence study for $h$-refinement for nonlinear waves in all combinations between $kh = \{1,3,6\}$ and $\varepsilon / \varepsilon_{\max} = \{10 \%, 50\%, 90 \%\}$.}
    \label{fig:h_convergence}
\end{figure}

The results of the $h$-refinement convergence study are presented in Figure~\ref{fig:h_convergence}, using basis functions of order $p = \{1, 2, 3, 4, 5 \}$. The parameter $h_{\max}$ denotes the maximum length of the element in the computational mesh. Due to the gradient recovery of $w$ from $\phi$, it is well known that an order of algebraic convergence is typically lost, \cite{engsigkarup2016stabilised}, leading to an expected convergence rate of $\mathcal{O}(h_{\max}^p)$. This expected behavior is clearly observed in the results. Furthermore, increasing the nonlinearity of the wave is shown to increase the complexity of the problem, requiring finer meshes to attain comparable levels of accuracy.

The $p$-convergence study is presented in Appendix A, in Figure \ref{fig:p_convergence}, where spectral, i.e., exponential, convergence is observed. The test is conducted on six different meshes with increasing resolution. Specifically, Mesh 1 has $h_{\max} = 0.25$, Mesh 2 has $h_{\max} = 0.125$, Mesh 3 has $h_{\max} = 0.0625$, etc. As observed previously, increasing nonlinearity results in a more challenging problem, requiring higher resolution to achieve comparable accuracy. \textbf{Remark:} Similar analyses have been performed for the LPF model, which exhibits improved convergence behavior, which is an expected outcome for a linearized model. These results are omitted here for the sake of brevity.

\subsubsection{Performance profiling}

It is well established that solving the Laplace problem constitutes the primary computational bottleneck in both LPF and FNPF wave models. Distributing the finite element domain across multiple cores inevitably affects performance. To this end, we profile the three most computationally intensive routines involved in time integration:  
1) \textit{LaplaceSolve}, which includes imposing boundary conditions, solving the linear system of equations, and performing gradient recovery (\textbf{Step 4.2, 4.3, 4.4, 4.5});  
2) \textit{EvaluateRHS}, responsible for evaluating the two free surface boundary conditions (\textbf{Step 4.6}); and  
3) \textit{LaplaceUpdate}, which updates the fluid domain and associated functions (\textbf{Step 4.1}).

Profiling is carried out on $N_c = \{1, \ldots, 300\}$ cores using two different meshes of polynomial order $p = 4$ with fixed vertical resolution $N_z = 4$, as detailed in Table \ref{tab:Profiling}. Mesh 2 is four times larger than Mesh 1, both in terms of the number of elements and degrees of freedom (DoF).

\begin{table}[H]
    \centering
    \caption{Performance profiling meshes with fixed $N_z=4$ and $p=4$.}
    \begin{tabular}{c r r r r} \hline
        Mesh & $N_{\Omega}$ & $N_{\Gamma^{\eta}}$ & $\text{DoF}(\mathcal{T}_{\Omega})$ & $\text{DoF}(\mathcal{T}_{\Gamma^{\eta}})$ \\ \hline
        1    & 32,768        & 8,192               & 1,114,112                            & 65,536                                   \\
        2    & 131,072       & 32,768              & 4,456,448                            & 262,144                                  \\ \hline
    \end{tabular}
    \label{tab:Profiling}
\end{table}

\begin{figure}[H]
    \centering
    \includegraphics[scale=0.95]{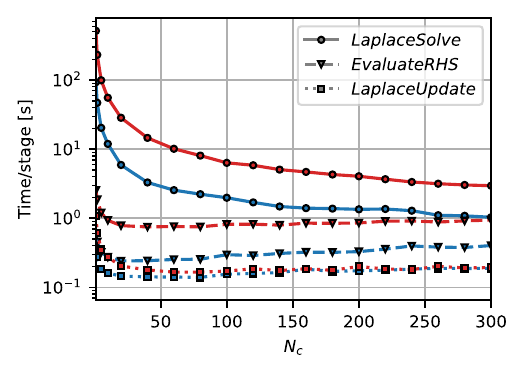}
    \caption{Timings for the three dominant routines in the FNPF model for Mesh 1 (blue) and Mesh 2 (red), using a polynomial order $p = 4$.}
    \label{fig:ProfilingTime}
\end{figure}

Figure \ref{fig:ProfilingTime} presents the core-averaged timing per time stage for the three routines. A substantial reduction in time per stage is observed; the increase is from $N_c = 1$ to approximately $N_c = 20$. Beyond this point, only the \textit{LaplaceSolve} routine continues to benefit from further parallelization, while \textit{EvaluateRHS} and \textit{LaplaceUpdate} begin exhibiting constant runtimes. Moving to a finer mesh (Mesh 2) increases the computational time for both \textit{LaplaceSolve} and \textit{EvaluateRHS}, as expected. In contrast, \textit{LaplaceUpdate} remains largely unaffected by mesh resolution, suggesting that its performance is dominated by communication overhead. The stagnation in \textit{EvaluateRHS} scalability is attributed to the relatively small size of the free surface mesh, which limits the parallel efficiency of this part. However, the overall solver remains scalable as long as \textit{LaplaceSolve} is the dominant cost factor. Moreover, from comparing the results of Mesh 1 and 2, it is clear that increasing the problem size allows for scalable behavior on a higher amount of cores due to \textit{LaplaceSolve} being the dominant cost factor for longer.

\begin{figure}[H]
    \centering
    \includegraphics[scale=0.95]{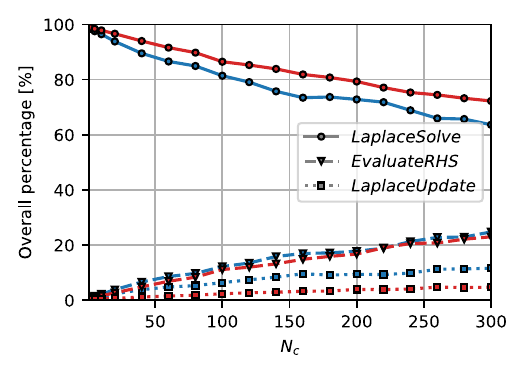}
    \caption{Percentage distribution of computational cost across the three dominant routines for Mesh 1 (blue) and Mesh 2 (red), using $p = 4$.}
    \label{fig:ProfilingPercent}
\end{figure}

Figure \ref{fig:ProfilingPercent} shows the relative percentage distribution of the total run time across the three routines. At low core counts (e.g., $N_c \leq 5$), \textit{LaplaceSolve} dominates the computational cost. However, as $N_c$ increases, its relative share decreases due to the constant cost of the remaining routines. For larger problems (as in Mesh 2), this imbalance is somewhat mitigated. Potential optimizations include reducing the core count assigned to the free surface evaluation or implementing more efficient iterative solvers for the nonlinear free surface problem. These optimizations are left for future work. \textbf{Remark:} The LPF formulation does not exhibit these performance limitations, as it avoids both the domain update step to solve \eqref{eq:Laplace_weak} and the nonlinear variational problem in \eqref{eq:free_surface_weak}, which reduces to a simple function evaluation in the linearized setting.

Based on this performance profiling of the code, it is confirmed that the Laplace solver constitutes the primary computational bottleneck. Consequently, the following strong- and weak-scaling analyses focus exclusively on this part of the code. Moreover, the convergent behavior showcased in the previous section has been confirmed on multicore hardware setups as well; however, is has been for conciseness.

\subsubsection{Parallel performance: Strong scaling}

To evaluate how run time performance is affected by increasing the number of cores for a fixed problem size (i.e., constant total workload), we perform strong scaling tests of the FNPF model. These tests assess both parallel efficiency and speedup as a function of the number of compute nodes, $N_n$, where each node comprises 20 cores. As a metric for strong scaling efficiency, we define $\gamma_s$ as
\begin{equation}\label{eq:StrongEfficiency}
    \gamma_s = \frac{T^b N_n^b}{T_{N_n} N_n},
\end{equation}
where $T^b$ and $N_n^b$ denote the runtime and node count of the baseline simulation, respectively, and $T_{N_n}$ is the runtime using $N_n$ nodes. In this study, the baseline corresponds to a single node configuration ($N_n^b = 1$). Strong scaling experiments are conducted using the three mesh configurations listed in Table \ref{tab:Strong}, with polynomial orders $p = \{3,4,5\}$ and a fixed vertical resolution of $N_z = 4$. The tests are executed on up to $N_n = 15$ nodes, equivalent to $N_c = 300$ cores, using a nonlinear periodic domain with $kh=2\pi$.

\begin{table}[H]
    \centering
    \caption{Strong scaling meshes with $N_z=4$.}
    \begin{tabular}{c r r c r r} \hline
        Mesh & $N_{\Omega}$ & $N_{\Gamma^{\eta}}$ & $p$  & $\text{DoF}(\mathcal{T}_{\Omega})$ & $\text{DoF}(\mathcal{T}_{\Gamma^{\eta}})$\\ \hline
            &  &  & 3           & 479,232 & 36,864 \\
        1    &  32,768  &  8,192  & 4                  & 1,114,112 & 65,536 \\
            &    & & 5                     & 2,150,400 & 102,400 \\ \hline
             &    &   & 3                  & 1,916,928 & 147,456 \\
        2    &  131,072  & 32,768  &  4                  & 4,456,448 & 262,144 \\
             &    &      & 5               & 8,601,600 & 409,600 \\\hline
             &    &   & 3                  & 3,833,037 & 294,849 \\
        3    & 262,088   & 65,522  &  4                  & 8,910,992 & 524,176 \\
             &    &  & 5                   & 17,199,525 & 819,025 \\\hline
    \end{tabular}
    \label{tab:Strong}
\end{table}

\begin{figure}[H]
    \centering
    \includegraphics[scale=0.95]{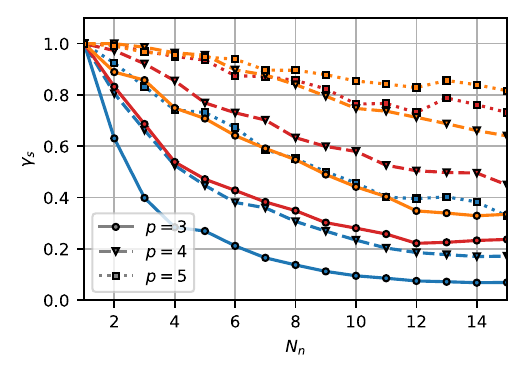}
    \caption{Strong scaling efficiency test. Mesh 1 (blue), Mesh 2 (red), and Mesh 3 (orange) with polynomial orders $p=\{3,4,5\}$.}
    \label{fig:StrongEfficiency}
\end{figure}

\begin{figure}[H]
    \centering
    \includegraphics[scale=0.95]{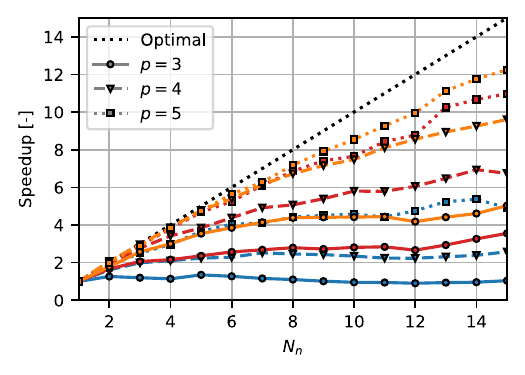}
    \caption{Speedup (compared to $N_n=1$ as baseline) for the strong scaling test. Mesh 1 (blue), Mesh 2 (red), and Mesh 3 (orange) with polynomial orders $p=\{3,4,5\}$.}
    \label{fig:StrongSpeedup}
\end{figure}

Figure \ref{fig:StrongEfficiency} shows the strong scaling efficiency using \eqref{eq:StrongEfficiency} for the three meshes. As expected, both increased mesh refinement and higher polynomial order improve the parallel efficiency of the model. This reflects better utilization of computational resources for larger problem sizes. The corresponding speedup curves are presented in Figure \ref{fig:StrongSpeedup}, alongside an ideal speedup line. The results indicate that the model approaches optimal speedup behavior as the problem size increases. For smaller problems, however, speedup plateaus early, revealing limitations in parallel efficiency due to insufficient workload per core.

\subsubsection{Parallel performance: Weak scaling}

To assess how the model scales with a constant workload per core, we performed a weak scaling analysis. The test domain is periodic, featuring nonlinear waves characterized by $kh = 2\pi$ and $L = 1$ [m]. The resolution is fixed in both the $y$- and $z$-directions with $N_y = N_z = 4$. The simulations are run on $N_n = \{1, \dots, 15\}$ nodes corresponding to $N_c = \{20, \dots, 300\}$ cores. To ensure a constant workload per core, the domain length is set as $L_x = N_c L$, which means that each core is responsible for simulating one wave within a domain of length $L$. The mesh configurations used for the tests are detailed in Table \ref{tab:Weak}, all using polynomial order $p = 5$. The largest simulation, executed on 300 cores, corresponds to a problem size of $\text{DoF}(\mathcal{T}_{\Omega}) = 80,640,000$. The weak scaling efficiency is defined as
\begin{equation}\label{eq:WeakEfficiency}
    \gamma_w = \frac{T^b}{T_{N_n}}.
\end{equation}

\begin{table}[H]
    \centering
    \caption{Weak scaling meshes with $N_y=N_z=4$ and polynomial order $p=5$. All reported quantities are per core, $N_c$.}
    \begin{tabular}{c r r c r r} \hline
        Mesh & $N_{\Omega}$ & $N_x$   & $\text{DoF}(\mathcal{T}_{\Omega})$ & $\text{DoF}(\mathcal{T}_{\Gamma^{\eta}})$\\ \hline
         1    &  128                & 4                          &  8,400                                  & 400  \\
         2    &  256                & 8                          &  16,800                                  & 800  \\
         3    &  512                & 16                          &  33,600                                  & 1,600  \\
         4    &  1,024                & 32                          &  67,200                                  & 3,200  \\
         5    &  2,048                & 64                          &  134,400                                  & 6,400  \\
         6    &  4,096                & 128                          &  268,800                                  & 12,800  \\
         \hline
    \end{tabular}
    \label{tab:Weak}
\end{table}

\begin{figure}[H]
    \centering
    \includegraphics[scale=0.95]{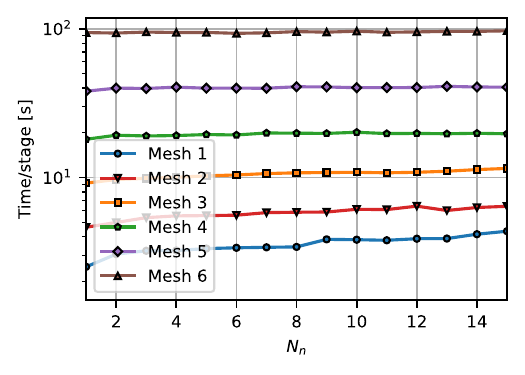}
    \caption{Weak scaling timings per stage on six meshes with polynomial order $p=5$.}
    \label{fig:WeakCPU}
\end{figure}

\begin{figure}[H]
    \centering
    \includegraphics[scale=0.95]{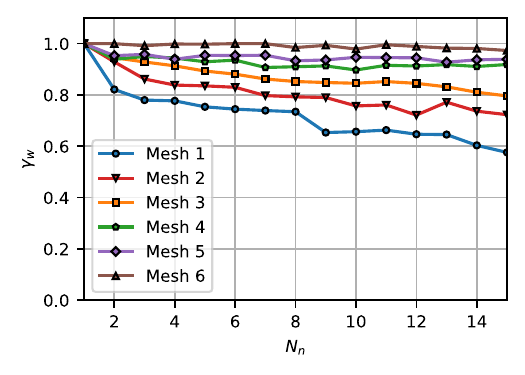}
    \caption{Weak scaling efficiency on six meshes with polynomial order $p=5$.}
    \label{fig:WeakEfficiency}
\end{figure}

Figure \ref{fig:WeakCPU} presents the run times for the various mesh configurations. As expected, increasing resolution leads to longer run times. However, more importantly, the runtime remains nearly constant as the number of cores increases, confirming the model's weak scaling behavior. This trend is further analyzed in Figure \ref{fig:WeakEfficiency}, which shows the weak scaling efficiency using \eqref{eq:WeakEfficiency}. The results indicate that higher mesh resolutions lead to improved parallel efficiency, as the computational workload per core becomes more substantial relative to communication overhead. \textbf{Remark:} Additional tests have been conducted for other polynomial orders, showing similar trends, as well as weak scaling tests for fixed meshes and increasing polynomial order. These results have been omitted for brevity.

\subsection{Model validation}

For all the validation cases, the simulations are run on a full node ($N_c = 20$ cores).

\subsubsection{Generation of high-order wave harmonics over a submerged bar}

First, we consider the classical bar test, which can be compared against the physical experimental laboratory data by Beji and Battjes, \cite{beji1994numerical}, denoted as BB94 in the following. Note that this is a wave flume test, thus two-dimensional; however, we extend the domain in the $y$-direction to make it fit our three-dimensional model. Mildly nonlinear waves are generated and propagated over a submerged bar, which generates high-order wave harmonics due to nonlinear shoaling effects (steepening and shortening of the waves). The domain is sketched out conceptually in Figure \ref{fig:FluidDomain}, where the specifications of the domain are: Generation zone ($0 \leq x \leq 8$). Target zone with the submerged bar ($8 \leq x \leq 30$). Absorption zone ($30 \leq x \leq 38$). The width of the domain is $0 \leq y \leq 1$. All units are [m]. The depth before and after the bar is $h = 0.4$ [m] and $h = 0.1$ [m] on top of the bar. The incline before and decline after the bar are 1:20 and 1:10, respectively. The input wave characteristics are set as $kh = 0.6725$ [-] and $\varepsilon/\varepsilon_{\max} \approx 7 \%$ given a wave period of $T = 2.018$ [s]. For the simulation, we used a free surface mesh of 7,600 elements, which is extruded with $N_z = 4$ elements in the vertical to give 22,800 prismatic elements for the representation of $\Omega$. Ultimately, this gives a problem size of 810,693 DoF and 62,361 DoF for $\Omega$ and $\Gamma^{\eta}$, respectively, when using basis functions of order $p=4$. The simulation is run for 25 wave periods. 

\begin{figure}[H]
    \centering
    \includegraphics[scale=0.95]{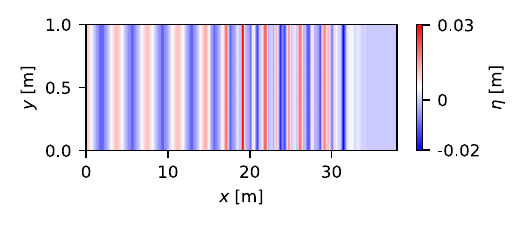}
    \caption{Contour plot of the free surface elevation at $t = 25T$ for the bar test.}
    \label{fig:BarTestFreeSurfaceElevation}
\end{figure}

\begin{figure}[H]
    \centering
    \includegraphics[scale=0.95]{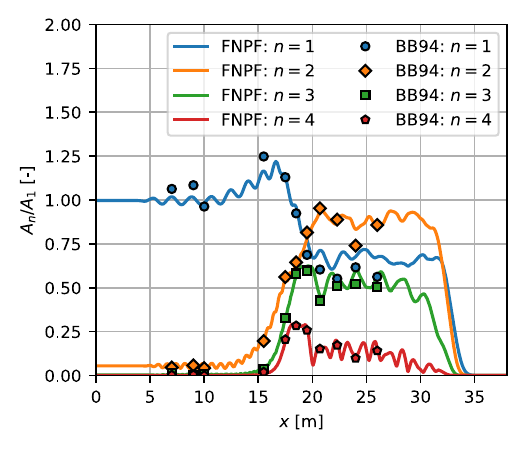}
    \caption{Harmonic analysis of the numerical simulation (FNPF) and experimental data (BB94) for the four most dominant harmonic modes in the bar test. $A_1 = H/2 = 0.01$ [m].}
    \label{fig:BarTestHarmonic}
\end{figure}

A contour plot of the free surface elevation, $\eta$, at the final time is shown in Figure \ref{fig:BarTestFreeSurfaceElevation}. Here, the high-order harmonic can be seen to take place when the wave reaches the submerged bar. In addition, a harmonic decomposition of the wave signal is carried out and presented in Figure \ref{fig:BarTestHarmonic}. The harmonic analysis is performed by fitting the time series, both experimental and numerical, using a least-squares approach to a sum of sinusoidal functions with argument $2\pi f_n t$, where $f_n = n/T$ and $n = \{1,\dots,4\}$. The figure illustrates how waves with dominant first-order components, along with minor second-order contributions, propagate toward the bar. Over and beyond the bar, second-order components begin to dominate the wave dynamics, with contributions from the first, third, and fourth harmonics. Finally, in the absorption zone, all wave energy is dissipated. The comparison between the numerical solution (FNPF) and the experimental data (BB94) shows good visual agreement.

\subsubsection{Linear and nonlinear interactions with a vertical cylinder}

Next, we investigate wave diffraction and the maximum free surface elevation, i.e., wave run-up, on a vertical surface-piercing cylinder subjected to both linear and nonlinear waves. The computational domain is symmetric about the $y=0$ plane, and the outer bounds are defined by $0 \leq x \leq 12L$ and $0 \leq y \leq 8L$. Wave generation and absorption zones are located in the regions $0 \leq x \leq 2L$ and $10L \leq x \leq 12L$, respectively. The vertical cylinder, of radius $R$, is centered at $(x, y) = (6L, 0)$. Each simulation is run for 25 wave periods.
The horizontal mesh is fully unstructured and graded to provide higher resolution near the cylinder wall boundary $\Gamma^{\text{cylinder}}$, with approximately eight elements per radius $R$. In contrast, far-field boundaries are meshed more coarsely, using roughly two elements per wavelength $L$. The vertical discretization uses $N_z = 4$ layers. The cylinder boundary is currently approximated using affine geometry, with a higher-order curvilinear representation left for future development. All relevant parameters for the domain, mesh, and wave conditions are summarized in Table \ref{tab:VerticalCylinder} for both the linear and nonlinear cases. In Appendix B, Figure \ref{fig:Cylinder_mesh}, the unstructured mesh used to obtain the linear solution can be seen.

\begin{table}[H]
    \centering
    \caption{Domain, mesh, and wave parameters for vertical cylinder simulations.}
    \begin{tabular}{l r r} \hline
         & Linear & Nonlinear \\ \hline
        Test case (TC) & 1 & 2 \\
        $R$ & $0.5$ [m]   & $0.1625$ [m] \\
        $L$ & $1$ [m]   & $2.73$ [m] \\
        $h$ & $1/(2\pi)$ [m]   & $0.45$ [m] \\
        $kh$ & $1$ [-]  & $1.036$ [-] \\
        $H$ & $0.01$ [m] & $0.053$ [m] \\
        $\varepsilon/\varepsilon_{\max}$ & $0\%$ & $19 \%$ \\
        $p$ & $4$ & $5$ \\
        $N_{\Omega}$ & 8,824 & 15,048 \\
        $N_{\Gamma^{\eta}}$ & 2,206 & 3,762 \\
        $\text{DoF}(\mathcal{T}_{\Omega})$ & 305,133 & 1,000,251 \\
        $\text{DoF}(\mathcal{T}_{\Gamma^{\eta}})$ & 17,949 & 47,631 \\  
        \hline
    \end{tabular}
    \label{tab:VerticalCylinder}
\end{table}

\begin{figure}[H]
    \centering
    \includegraphics[scale=0.90]{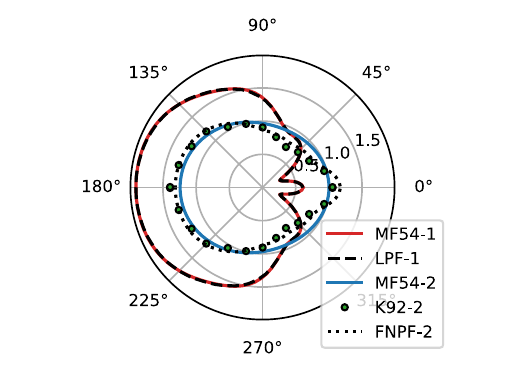}
    \caption{Maximum non-dimensional free surface elevation on the cylinder, $\eta_m/(2H)$, for the two simulations comparing numerical, experimental, and analytical solutions.}
    \label{fig:VerticalCylinder_solution}
\end{figure}

The results from both simulations are shown in Figure \ref{fig:VerticalCylinder_solution}, presented as a non-dimensional maximum free surface elevation on the cylinder boundary, $\eta_m/(2H)$. Here, $\eta_m$ is defined as
\begin{equation*}
    \eta_m(x_i, y_i) = \max_t \eta(x_i, y_i, t), \quad (x_i, y_i) \in \Gamma^{\text{cylinder}}.
\end{equation*}
For the linear test case (TC 1), the numerical results, denoted LPF-1, are compared with the analytical predictions of linear diffraction theory by \cite{maccamy1954wave}, denoted MF54-1. For the nonlinear test case (TC 2), the numerical results of the FNPF simulation, denoted FNPF-2, are compared with the digitized experimental data of \cite{kriebel1992nonlinear}, denoted K92-2, which were digitized for comparison. We also include equivalent linear results for this configuration, denoted MF54-2.

Figure \ref{fig:VerticalCylinder_solution} shows good agreement between the numerical and analytical results for the linear case (LPF-1 vs. MF54-1), particularly in the region downstream of the cylinder. For the nonlinear test case, the numerical solution (FNPF-2) aligns well with the experimental data (K92-2), consistent with previous results reported in the literature, e.g., see \cite{ducrozet2010high}. Notably, nonlinear effects (compared to MF54-2) manifest as increased run-up at $0^\circ$ and $180^\circ$, and a reduced run-up at intermediate angles, around $45^\circ$ and $315^\circ$.

\subsubsection{Wave interactions with a V-shaped breakwater}

As a final validation test case, we consider wave interactions with a breakwater structure, following the setup described in \cite{chang2012accurate}, for incident waves approaching at angles of $\theta = \{0^\circ, 30^\circ, 60^\circ\}$. The computational domain is defined as follows: From the origin, $(x,y) = (0,0)$, two breakwater arms extend a length of $b = 76.2$ [m] and form an internal angle of $2\beta = \pi/3$ [rad] ($60^\circ$). We refer to the upper and lower arms as the upper breakwater arm (UBA) and lower breakwater arm (LBA). At the ends of the arms, inward-facing curved end-walls are placed, each spanning an angle of $2\gamma = 2\pi/9$ [rad] ($40^\circ$). While the original study by \cite{chang2012accurate} treated the breakwaters as infinitely thin, we adopt a more physically realistic representation by assigning them a finite thickness of $l = 1$ [m].

\begin{figure}[H]
    \centering
    \includegraphics[scale=0.95]{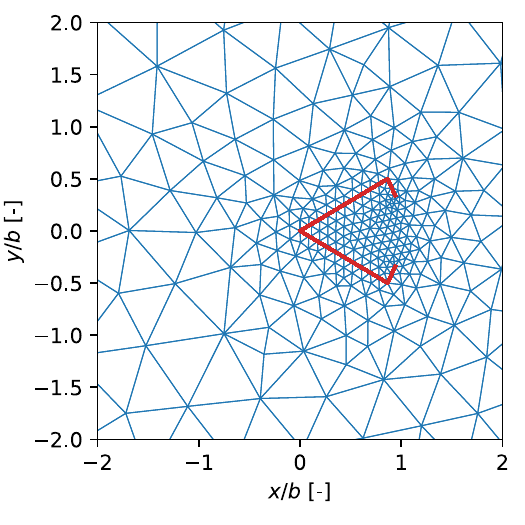}
    \caption{Free surface mesh for $\theta = 0^\circ$ around the breakwater.}
    \label{fig:Breakwater_mesh}
\end{figure}

The computational domain spans the region $-L_2 - 2b \leq (x,y) \leq 2b + L_2$, completely enclosing the target region defined as $-2b \leq x \leq 2b$ and $-L_2 - 2b \leq y \leq 2b + L_2$. Here, $L_2 = 260$ [m] $ \approx 2L$, with $L$ denoting the wavelength of the incoming wave. A generation zone is defined for $-L_2 - 2b \leq x \leq -2b$, and an absorption zone for $2b \leq x \leq 2b + L_2$. The water depth is set to $h = 13.4$ [m] throughout. For the spatial discretization, a fully unstructured horizontal mesh is employed. The mesh is coarsened toward the domain ends (approximately one element per wavelength, $L$) and refined near the breakwaters (approximately ten elements per breakwater arm length, $b$). This results in 754 triangular elements on the free surface, $\mathcal{T}_{\Gamma^\eta}$, which are extruded into four vertical layers, yielding a total of 3,016 prismatic elements in the fluid domain, $\mathcal{T}_{\Omega}$. 
Using a polynomial order of $p=5$, the resulting discrete system contains 202,020 degrees of freedom (DoF). The free surface mesh can be seen in Figure \ref{fig:Breakwater_mesh}.

\begin{figure}[H]
    \centering
    \includegraphics[scale=0.95]{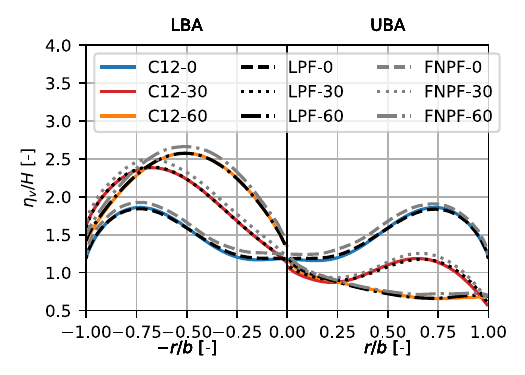}
    \caption{Breakwater solution (non-dimensional free surface variation) on the UBA and LBA for $\theta = \{0^\circ, 30^\circ, 60^\circ\}$ using proposed linear (LPF-$\theta$) and nonlinear (FNPF-$\theta$) numerical model and results from \cite{chang2012accurate} (C12-$\theta$).}
    \label{fig:Breakwater_solution_arms}
\end{figure}

\begin{figure}[H]
    \centering
    \includegraphics[scale=0.95]{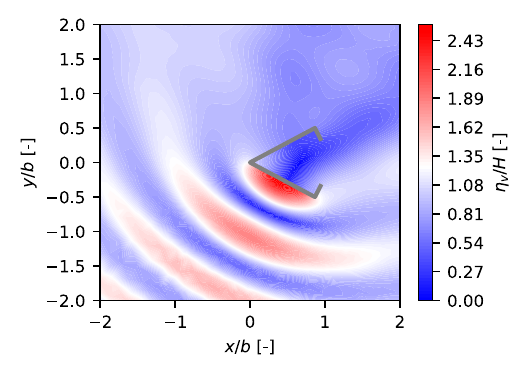}
    \caption{Breakwater solution (normalized free surface variation) at $\theta = 60^\circ$ angle of attack.}
    \label{fig:Breakwater_solution_60}
\end{figure}

The wave input parameters are taken from \cite{chang2012accurate}, with a non-dimensional wavenumber $kh = 0.65$ [-], corresponding to $L = 129.5$ [m]. The linear -- and arbitrary -- wave height is set to $H = 1.295$ [m], and the simulations are run for 10 wave periods. For each incident angle, $\theta = \{0^\circ, 30^\circ, 60^\circ\}$, we compute the time-independent free surface elevation variation, defined at each grid point $(x_i, y_i)$ as
\begin{equation*}
   \eta_v =  \eta_v(x_i, y_i) = \max_t \eta(x_i, y_i, t) - \min_t \eta(x_i, y_i, t).
\end{equation*}
The corresponding non-dimensional quantity $\eta_v/H$ is visualized as contour plots in Appendix C (Figures \ref{fig:Breakwater_solution_0}, \ref{fig:Breakwater_solution_30}, and \ref{fig:Breakwater_solution_60}). In Figure \ref{fig:Breakwater_solution_arms}, we present $\eta_v/H$ along the UBA and LBA for the three wave angles, comparing our numerical results (LPF-$\theta$) with digitized reference data from \cite{chang2012accurate} (C12-$\theta$). Despite small differences arising from the use of finite-thickness structures, the discretizations of the domain, and minor wave reflections at the boundaries, the numerical solutions show good qualitative agreement with the reference results. This highlights that the incident angle impacts the wave elevation in the vicinity of the breakwater structure, e.g., which may be of interest in sea state or load estimations on the structure. Added to the figure, we also show the results from the FNPF model (FNPF-$\theta$) with the same wave input parameters; however, with an incident wave of $\varepsilon/\varepsilon_{\max} = 18\%$. From the figure, the added nonlinearity is visually evident in the maximum variance of the free surface elevation on the structure, especially on the LBA side. For $\theta = 60^\circ$ on the UBA side, only a minor change can be observed; however, that is also the side completely sheltered from the incident waves.

\section{Conclusion}\label{sec:con}

In this paper, we presented a new parallel implementation of a linear and fully nonlinear potential flow model, FNPF-SEM, to simulate wave propagation and wave-structure interaction. The model employs high-order (spectral) finite elements within the open-source Firedrake framework, leveraging its native MPI-based parallelism to handle large-scale computations efficiently.
We detailed the mathematical formulation, including the governing equations, domain definitions, and numerical discretization. Parallel performance was assessed through code performance profiling for strong- and weak-scaling studies, which demonstrated good parallel efficiency, particularly for large-scale simulations. Profiling confirmed that the Laplace problem constitutes the dominant computational cost.
The numerical accuracy of the implementation was verified through convergence studies, confirming the expected algebraic and spectral convergence under mesh refinement and increased polynomial order, respectively. Finally, the model was validated against analytical solutions and experimental data for a range of wave propagation and wave-structure interaction scenarios, including high-order harmonic generation over a submerged bar and linear and nonlinear wave interaction with a vertical cylinder and a V-shaped breakwater. These results confirm the accuracy and applicability of the proposed framework. In ongoing work, the tool will be expanded to handle regional-scale wave propagation for accurate sea state estimation and consider areas of irregular coastlines, exploiting the advantage of unstructured meshes.

\section*{Acknowledgments}
The research was carried out at the Technical University of Denmark (DTU) in the Department of Applied Mathematics and Computer Science with computational resources provided by the DTU Computing Center (DCC). Lastly, the authors express their gratitude to the Firedrake team led by Professor David Ham for their support.

\section*{Declarations}
The authors declared that they had no potential conflicts of interest with respect to the research, authorship, and/or publication of this article.

\section*{Funding}
The authors disclosed receipt of the following financial support for the research, authorship, and/or publication of this article: This work partly contributes to the activities of the research project: ‘‘A new digital twin concept for floating offshore structures’’ supported by COWIfonden (Grant no. A-165.19) and partly the PhD-project of JV: ‘‘New Advanced Simulation Techniques for Wave Energy Converts’’.





\bibliographystyle{SageH}

\bibliography{references.bib}  





\newpage
\section{Appendix}\label{sec:app}

\subsection{Appendix A: \textit{p}-convergence study}

\begin{figure}[H]
    \centering
    \includegraphics[scale=0.90]{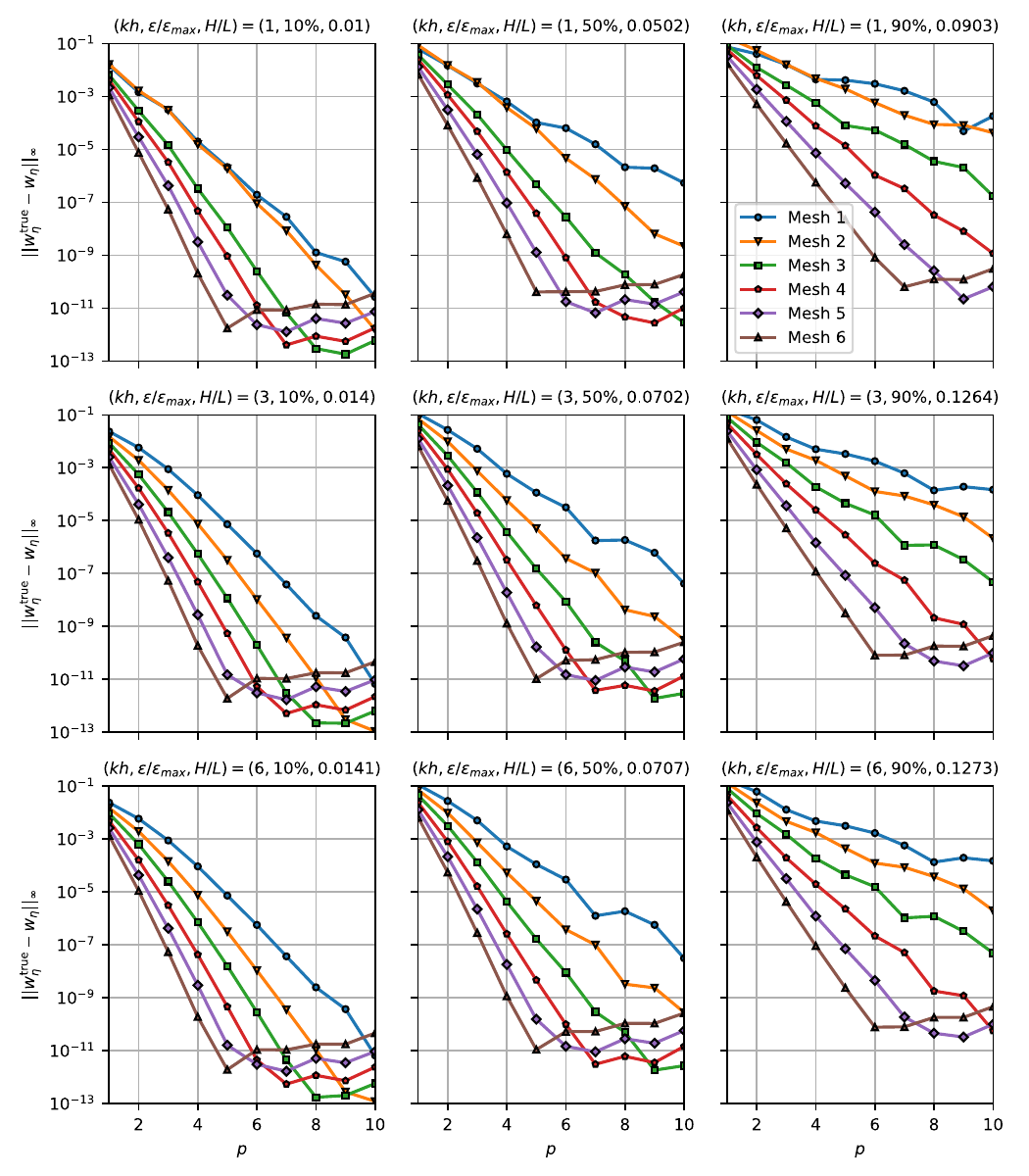}
    \caption{Convergence study for $p$ for nonlinear waves in all combinations between $kh = \{1,3,6\}$ and $\varepsilon / \varepsilon_{\max} = \{10 \%, 50\%, 90 \%\}$.}
    \label{fig:p_convergence}
\end{figure}

\newpage
\subsection{Appendix B: Vertical cylinder mesh}
\begin{figure}[H]
    \centering
    \includegraphics[scale=0.95]{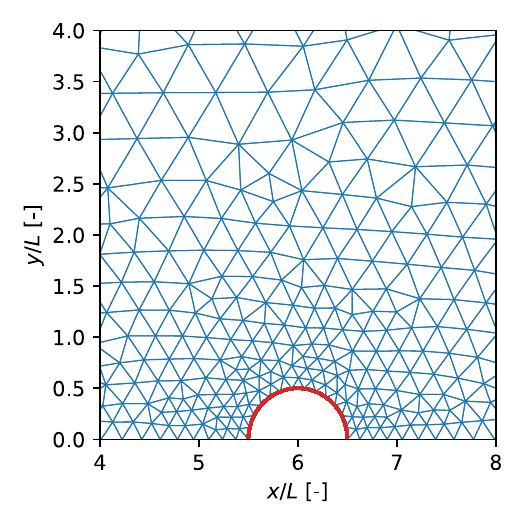}
    \caption{Mesh used for the linear simulation as described in Table \ref{tab:VerticalCylinder}. Zoomed in around the cylinder.}
    \label{fig:Cylinder_mesh}
\end{figure}

\newpage
\subsection{Appendix C: Breakwater solutions}

\begin{figure}[H]
    \centering
    \includegraphics[scale=0.95]{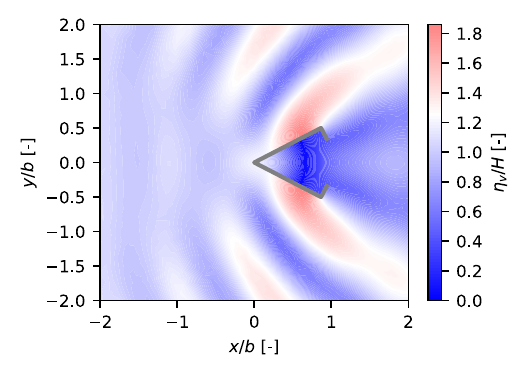}
    \caption{Breakwater solution (normalized free surface variation) at $\theta =0^\circ$ angle of attack.}
    \label{fig:Breakwater_solution_0}
\end{figure}

\begin{figure}[H]
    \centering
    \includegraphics[scale=0.95]{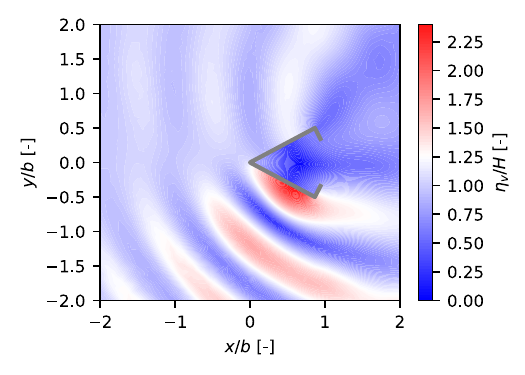}
    \caption{Breakwater solution (normalized free surface variation) at $\theta =30^\circ$ angle of attack.}
    \label{fig:Breakwater_solution_30}
\end{figure}

\end{document}